\documentclass{article}
\usepackage{amsmath}
\usepackage{amsthm}
\usepackage{amstext}
\usepackage{amsopn}
\usepackage{graphicx}

\usepackage{amsfonts}
\usepackage{amssymb}
\usepackage{mathrsfs}
\theoremstyle{theorem}
\newtheorem{theorem}{Theorem}

\theoremstyle{definition}

\newtheorem{example}[theorem]{Example}

\begin{document}
\date{April 5, 2017}
\title{A note on Dodgson's determinant condensation algorithm}
%\author{} %%leave blank in initial submission to allow for double blind reviewing
\author{Hou-biao Li\footnote{Corresponding Author, Email: lihoubiao0189@163.com}, Hong Li, Ting-zhu Huang\\
{\small School of Mathematical Sciences, University of Electronic Science}\\
{\small and Technology of China, Chengdu, 611731, P.R. China.}\\
{\small and Technology of China, Chengdu, P.R. China}}

\maketitle
%Notes authors should leave bios out altogether.

%\begin{biog} %comment out in initial submission to allow for double blind reviewing
%\item[\biogpic{\includegraphics[width=84pt]{Woodrow_Wilson.pdf}}Woodrow Wilson] (twwilson@princeton.edu) received his PhD in history and political science from Johns Hopkins University. He held visiting positions at Cornell and Wesleyan before joining the faculty at Princeton, where he was eventually appointed president of the university.  Among his proudest accomplishments was the abolition of eating clubs at Princeton on the grounds that they were elitist.

%\item[\biogpic{\includegraphics[width=84pt]{Herbert_Hoover.pdf}}Herbert Hoover] (hchoover@stanford.edu) entered Stanford University in 1891, after failing all of the entrance exams except mathematics.  He received his BS degree in geology in 1895, spent time as a mining engineer, then was appointed by his co-author to the U.S. Food Administration and the Supreme Economic Council, where he orchestrated the greatest famine relief efforts of all time.
%\end{biog}
\begin{abstract}
Recently, the Dodgson's determinant condensation algorithm was revisited in many papers [\emph{College Math. Journal 42(1)(2011): 43--54, College Math. Journal 38(2)(2007): 85--95}, \emph{Math Horizons 14(2)(2006): 12--15},etc.]. This method is attractive, but there also exist some shortcomings. In this paper, a symbolic algorithm and the corresponding problems are discussed to overcome these shortcomings. Numerical experiments show that this symbolic modified algorithm has highly sensitivity on initial configuration of the matrix in condensation process.
\end{abstract}

\section{Introduction}
In this paper, we mainly consider an $n\times n$ matrix $A=(a_{ij})$ with elements $a_{ij}$ real or complex number. Its determinant is denoted by $|A|$, also written $\det A$. For the calculation of determinants, the Dodgson's determinant condensation algorithm was recently revisited in many papers \cite{LAA2008, C2011, C2006, C2007}. This method is uncomplicated by the calculation of $2\times 2$ determinants and unique utilization of matrix condensation techniques which has promise for parallel computing architectures. However, there does not seem to be sensitivity on initial configuration of the matrix in condensation process for intermediate matrices, which effects its application.

By observing the computation process of the Dodgson's determinant condensation algorithm, a modified symbolic algorithm with highly sensitivity to initial conditions is presented in this paper and the corresponding problems are discussed to overcome the above shortcoming.

\section{The Dodgson's condensation algorithm}

As it is well-known, Dodgson's condensation algorithm is also referred to as the Lewis' method (see \cite{LAA2008,Do1866,C2011,C2006,C2007}).
This algorithm is one of considerable computational simplicity, achieved by restricting itself entirely to the calculation of $2\times 2$ determinants, which usually consists of the following steps \cite{C2006}:
\begin{enumerate}
  \item Use elementary row and column operations to remove all zeros from the interior of $A$.
  Here, the interior of an $n\times n$ ($n>3$) matrix $A$, or intA, is the $(n-2)\times(n-2)$ consecutive minor that results
  when the first and last rows and columns of matrix $A$ are deleted.
  \item Find the $2\times 2$ determinant for every four adjacent terms to form a new $(n-1)\times(n-1)$ matrix $A^{(n-1)}$.
  \item Repeat this step to produce an $(n-2)\times (n-2)$ matrix, and then divide each term by the corresponding entry in the interior of
        the original matrix $A$, to obtain matrix $A^{(n-2)}$.
  \item Continue ``condensing" the matrix down, until a single number $A^{(1)}$ is obtained. This final number will be $\det A$.
\end{enumerate}

Obviously, each iterate $A^{(k)}$ is a condensation of the previous iterate $A^{(k+1)}$. In fact, we have that

\begin{theorem}\label{Th2.2}(Dodgson's Condensation Theorem, \cite{C2007}).
Let $A$ be an $n\times n$ matrix. After $k$ successful condensations, Dodgson's method produces the matrix
\begin{equation}\label{2.1}
{A^{(n - k)}} = \left[ {\begin{array}{*{20}{c}}
{|{A_{1 \ldots k + 1,1 \ldots k + 1}}|}&{|{A_{1 \ldots k + 1,2 \ldots k + 2}}|}& \cdots &{|{A_{1 \ldots k + 1,n - k \ldots n}}|}\\
{|{A_{2 \ldots k + 2,1 \ldots k + 1}}|}&{|{A_{2 \ldots k + 2,2 \ldots k + 2}}|}& \cdots &{|{A_{2 \ldots k + 2,n - k \ldots n}}|}\\
 \vdots & \vdots & \ddots & \vdots \\
{|{A_{n - k \ldots n,1 \ldots k + 1}}|}&{|{A_{n - k \ldots n,2 \ldots k + 2}}|}& \cdots &{|{A_{n - k \ldots n,n - k \ldots n}}|}
\end{array}} \right],
\end{equation}
whose entries are the determinants of all $(k+1)\times(k+1)$ submatrices of $A$, where $A_{i\ldots j,k\ldots l}$ denotes the
submatrix composed of rows $i,i+1,\ldots,j$ and columns $k,k+1,\ldots, l$ of $A$.
\end{theorem}

If Dodgson's method terminates successfully, it computes the determinant of an $n\times n$ matrix using
\begin{equation}\label{1.5}
2\left[\left(n-1\right)^{2}+\left(n-2\right)^{2}+\cdots+1^{2}\right]=\frac{2}{3}{n^3}-{n^2} + \frac{1}{3}n\approx \frac{2}{3}n^{3}
\end{equation}
multiplications, which has the same multiplications as the Chi\`{o}'s condensation method \cite{F75}. However, division presents Dodgson's method with a huge drawback, if the central element of $A^{(k)}$ for any $k>2$ is null, the reduction method cannot be directly
applied. Fortunately, we sometimes may, for instance, before applying the method, permutate cyclically the lines (or columns) of
$A$ in such a way that the new central element turns to be a non-null element. But, as mentioned in \cite{C2011}, there also exist some matrices
(see the following example \ref{E1.1}) making this technology noneffective.

\begin{example}\label{E1.1} (\cite{C2011}). No combination of row or column swaps allows Dodgson's method to
compute the determinant of
\begin{equation}\label{1.6}
A = \left[ {\begin{array}{*{20}{c}}
1&0&3&0\\
0&{ - 1}&0&1\\
1&1&2&0\\
0&2&0&1
\end{array}} \right],
\end{equation}
because there will always be a zero in the interior of $A$.
\end{example}

\section {Methods for overcoming the Dodgson's drawback}

In this section, we will review some new methods to overcome the above Dodgson's problem caused by zeros in the interior
of a matrix in Section 1.

%\subsection{Method one---zero elements symbolic algorithm}

Firstly, D. Leggett {\em et al.} presented a double-crossing method to solve this problem in the literature \cite{C2011} whenever the non-zero element is immediately above, below, left, right, or catty-corner to the zero; that is, the non-zero element is adjacent to the zero. But if the zero appears in a $3\times 3$ block of zeros, then the double-crossing method may fail and its corresponding calculations are very complex.

Since a determinant is a continuous function on its elements, one may choose some symbolic variables to replace the zeros which occur in the interior of the matrix $A^{(n-k)}$ during the computing process and then compute $A^{(n-(k+1))}$ and $A^{(n-(k+2))}$ as usual.
Finally, let these symbolic variables approach to zero in the matrix $A^{(2)}$ or $A^{(1)}$, and obtain the determinant value.
This technology was even used to find the inverse of a tridiagonal matrix in \cite{E2006,Li2010}.
To illustrate this method, let us firstly consider the following fourth order determinant.

\begin{example}\label{E2.1}. Let $A^{(4)}$ be the matrix $A$ of Example \ref{E1.1}. We have one zero
element in the interior, at $(i,j)=(2,3)$. As set forth in the Double-Crossing Theorem \cite{C2011},
this affects element $(1,2)$ of $A^{(2)}$, i.e.,
\[
{A^{(3)}} = \left[ {\begin{array}{*{20}{c}}
{ - 1}&3&3\\
1&{ - 2}&{ - 2}\\
2&{ - 4}&2
\end{array}} \right] \to {A^{(2)}} = \left[ {\begin{array}{*{20}{c}}
1&?\\
0&{ - 6}
\end{array}} \right].
\]

Next, let us solve this problem by the above symbolic variables method.
Let the zero element of $A^{(4)}$, at $(i,j)=(2,3)$, be denoted symbolically by $x$ and continue computing as usual.

\begin{equation}\label{2.4}
{A^{(4)}} = \left[ {\begin{array}{*{20}{c}}
1&0&3&0\\
0&{ - 1}&x&1\\
1&1&2&0\\
0&2&0&1
\end{array}} \right] \to {A^{(3)}} = \left[ {\begin{array}{*{20}{c}}
{ - 1}&3&3\\
1&{ -(x+2)}&{ - 2}\\
2&{ - 4}&2
\end{array}} \right]
 \to {A^{(2)}} = \left[ {\begin{array}{*{20}{c}}
{1 - x}&3\\
{2x}&{ - x - 6}
\end{array}} \right];
\end{equation}

\begin{equation}\label{2.5}
{A^{(1)}} = \left( {\frac{{\left| {{A^{(2)}}} \right|}}{{ - (x + 2)}}} \right) = \frac{{{x^2} - x - 6}}{{ - (x + 2)}} = 3 - x.
\end{equation}
In \eqref{2.4} and \eqref{2.5}, taking the limit as $x$ approaches zero, we obtain the determinant of $A^{(4)}$.
\end{example}

Comparing the above process with the double-crossing method, it is obvious that the symbolic method seems more simple, since
it does not need new sub-matrices or determinant calculates.
Specially, it preserves the general simplicity and the idea of Dodgson's method.

However, the above method does not yield the correct result in some cases if we replace the zero entries
in the interior of intermediate matrix $A^{(n-k)}$ with symbolic variables instead of changing the corresponding entries of the original matrix. For example,
\begin{example}\label{E2.2}(\cite{CMJ11}). Consider the following matrix
\begin{equation}\label{2.6}
A = \left[ {\begin{array}{*{20}{c}}
3&{ - 2}&1&2\\
{ - 1}&4&4&1\\
3&3&3&4\\
2&5&2&{ - 1}
\end{array}} \right].
\end{equation}
The original Dodgson's method yields
\begin{equation}\label{2.7}
{A^{(3)}} = \left[ {\begin{array}{*{20}{c}}
{10}&{ - 12}&{ - 7}\\
{ - 15}&0&{13}\\
9&{ - 9}&{ - 11}
\end{array}} \right],
\end{equation}
whose interior element, $a_{22}$, is zero. Suppose we replace $a_{22}$ with the symbolic variable $x$ to get the matrix $B^{(3)}$, i.e.,
\begin{equation}\label{2.8}
{B^{(3)}} = \left[ {\begin{array}{*{20}{c}}
{10}&{ - 12}&{ - 7}\\
{ - 15}&x&{13}\\
9&{ - 9}&{ - 11}
\end{array}} \right].
\end{equation}
And then we obtain the $B^{(2)}$ by applying Dodgson's condensation to $B^{(3)}$:
\[
{B^{(2)}} = \left[ {\begin{array}{*{20}{c}}
{\frac{{5x - 90}}{2}}&{\frac{{7x - 156}}{4}}\\
{45 - 3x}&{\frac{{117 - 11x}}{3}}
\end{array}} \right] \Rightarrow {B^{(1)}} = {{\left( {\frac{{267}}{4}x - \frac{{47}}{{12}}{x^2}} \right)} \mathord{\left/
 {\vphantom {{\left( {\frac{{267}}{4}x - \frac{{47}}{{12}}{x^2}} \right)} x}} \right.
 \kern-\nulldelimiterspace} x} = \frac{{267}}{4} - \frac{{47}}{{12}}x.
 \]
Taking the limit as $x$ approaches zero, we have
\[
\mathop {\lim }\limits_{x \to 0} \det (B) = \mathop {\lim }\limits_{x \to 0} (\frac{{267}}{4} - \frac{{47}}{{12}}x) = \frac{{267}}{4} \ne 213 = \det A.
\]
\end{example}

A very natural question is what happens for the above process. In fact, the computation of determinants differs from the literature \cite{E2006,Li2010}. Because the inversion of a
general nonsingular matrix is unique, however, the $n$-order determinants with the same value are ubiquitous. For example, the following matrices have
the same intermediate matrix \eqref{2.7} in the computing process of the Dodgson's method.
\[
{A_1} = \left[ {\begin{array}{*{20}{c}}
1&{ - 2}&1&2\\
3&4&4&1\\
6&3&3&4\\
7&5&2&{ - 1}
\end{array}} \right],\;\;\;{A_2} = \left[ {\begin{array}{*{20}{c}}
{0.25}&9&{12}&{4.75}\\
{ - 1}&4&4&1\\
3&3&3&4\\
2&5&2&{ - 1}
\end{array}} \right],
\]

\[{A_3} = \left[ {\begin{array}{*{20}{c}}
3&{ - 2}&1&{0.25}\\
{ - 1}&4&4&{ - 6}\\
3&3&3&{ - 1.25}\\
2&5&2&{ - 4.5}
\end{array}} \right],
{A_4} = \left[ {\begin{array}{*{20}{c}}
{ - 45}&{ - 5}&{ - 2}&{0.5}\\
{ - 7}&{ - 1}&2&3\\
{ - 8}&1&{ - 2}&{3.5}\\
{ - 25}&2&{ - 13}&{28.25}
\end{array}} \right],
\]

\[{A_5} = \left[ {\begin{array}{*{20}{c}}
0&1&5&{ - 8.5}\\
{ - 10}&2&{ - 2}&2\\
{ - 7.5}&3&{ - 3}&{ - 3.5}\\
{ - 8}&2&{ - 5}&{\frac{{ - 13}}{6}}
\end{array}} \right], {A_6} = \left[ {\begin{array}{*{20}{c}}
0&1&4&{\frac{{31}}{{12}}}\\
{ - 10}&6&{12}&6\\
{\frac{{35}}{{12}}}&{ - \frac{1}{4}}&{ - \frac{1}{2}}&{\frac{5}{6}}\\
{ - 34}&6&{48}&{ - 58}
\end{array}}\right].
\]
However, their determinants are not the same:
\[\begin{array}{l}
\det {A_1} = \det {A_2} = \det {A_3} = 213,\\
\det {A_4} = 5665.5,\\
\det {A_5} = 451.5,\\
\det {A_6} = 2073,
\end{array}\]
which show that there doesn't seem to be sensitivity on initial configuration of the matrix in the Dodgson's determinant condensation calculations, see Figure 1.

\begin{figure}[htbp]\label{figure1}
\centering
\includegraphics[scale=0.75]{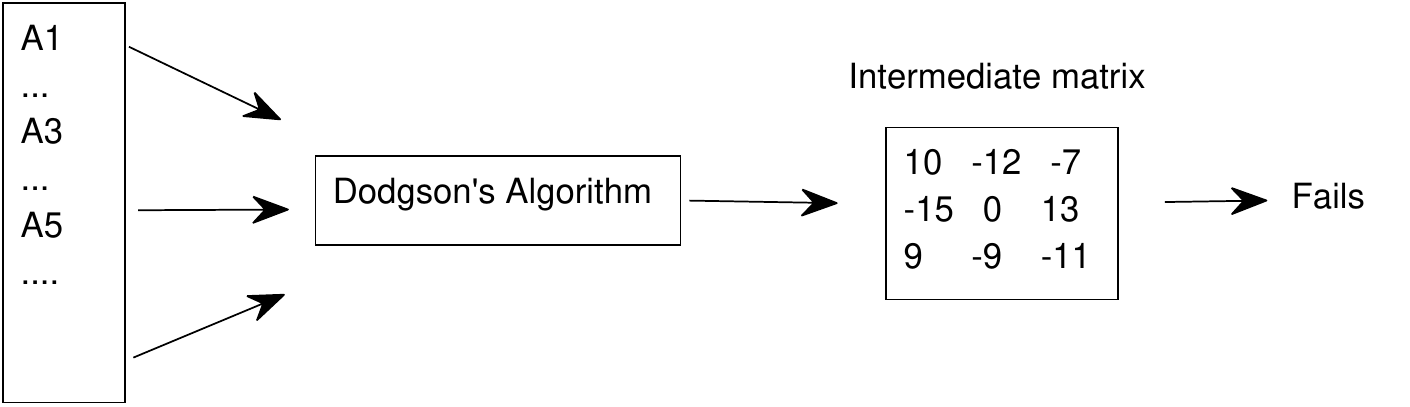}
\caption{The non-sensitivity of the Dodgson＊s algorithm to ``initial value".}
\end{figure}

\section{Modified Dodgson's condensation Algorithm}

To solve the above problem, there exist several tactics.

(1) Firstly, we may add a symbolic variable $x$ to an entry of the $2\times 2$ submatrix of $\textbf{original}$ matrices $A_i$ ($i=1,\ldots 6$)
whose determinant leads to the zero entry of $A^{(3)}$ (see \eqref{2.7}). Continuing with Dodgson's method, we will have the correct results,
as shown below. For example, for $A_1$ and $A_4$, we have respectively that

(i) \[\begin{array}{l}
{{\bar A}_1} = \left[ {\begin{array}{*{20}{c}}
1&{ - 2}&1&2\\
3&4&4&1\\
6&3&{3 + x}&4\\
7&5&2&{ - 1}
\end{array}} \right] \Rightarrow \bar A_1^{(3)} = \left[ {\begin{array}{*{20}{c}}
{10}&{ - 12}&{ - 7}\\
{ - 15}&{4x}&{13 - x}\\
9&{ - 9 - 5x}&{ - 11 - x}
\end{array}} \right]\\
 \Rightarrow \bar A_1^{(2)} = \left[ {\begin{array}{*{20}{c}}
{10x - 45}&{10x - 39}\\
{45 + 13x}&{39 - 9x}
\end{array}} \right] \Rightarrow \bar A_1^{(1)} = \frac{{852x - 220{x^2}}}{{4x}} = 213 - 55x;
\end{array}\]

(ii)\[\begin{array}{l}
{{\bar A}_4} = \left[ {\begin{array}{*{20}{c}}
{ - 45}&{ - 5}&{ - 2}&{0.5}\\
{ - 7}&{ - 1 + x}&2&3\\
{ - 8}&1&{ - 2}&{3.5}\\
{ - 25}&2&{ - 13}&{28.25}
\end{array}} \right] \Rightarrow \bar A_4^{(3)} = \left[ {\begin{array}{*{20}{c}}
{10 - 45x}&{ - 12 + 2x}&{ - 7}\\
{ - 15 + 8x}&{ - 2x}&{13}\\
9&{ - 9}&{ - 11}
\end{array}} \right]\\
 \Rightarrow \bar A_4^{(2)} = \left[ {\begin{array}{*{20}{c}}
{74x + 180}&{6x - 78}\\
{135 - 54x}&{ - \frac{{117}}{2} - 11x}
\end{array}} \right] \Rightarrow \bar A_4^{(1)} = \frac{{ - 11331x - 490{x^2}}}{{ - 2x}} = \frac{{11331}}{2} + 245x.
\end{array}\]

Taking the limit as $x$ approaches zero, we have
\[\begin{array}{l}
\mathop {\lim }\limits_{x \to 0} \det ({{\bar A}_1}) = \mathop {\lim }\limits_{x \to 0} (213 - 55x) = 213;\\
\mathop {\lim }\limits_{x \to 0} \det ({{\bar A}_4}) = \mathop {\lim }\limits_{x \to 0} (\frac{{11331}}{2} + 245x) = 5665.5.
\end{array}\]
That is to say, the symbolic modified Dodgson＊s algorithm is sensitivity on ``initial value" in the condensation calculation process, see Figure 2.

\begin{figure}[htbp]\label{figure:2}
\centering
\includegraphics[scale=0.75]{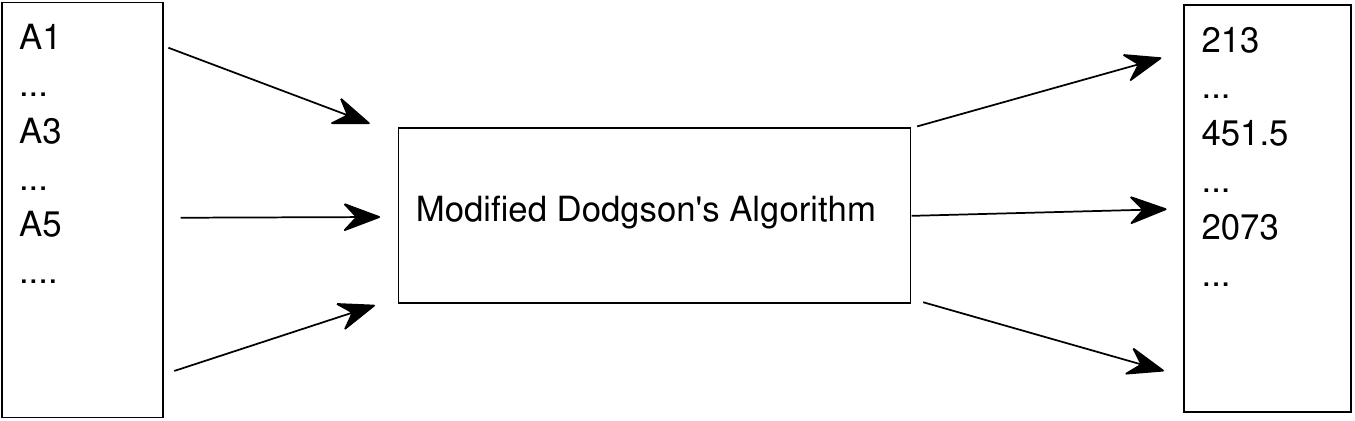}
\caption{The sensitivity of the symbolic modified Dodgson＊s algorithm to ``initial value".}
\end{figure}

Therefore, the sensitivity of the algorithm to ``initial value" is very important in algorithm design. Next, let us look at the following specific example.
\[{A_7} = \left[ {\begin{array}{*{20}{c}}
1&{ - ai}&1&2\\
3&3&{a + bi}&1\\
6&3&{a + i}&4\\
7&5&2&{ - 1}
\end{array}} \right],\]
where the $a$ changes from $2$ to $-3$, the $b$ changes from $0$ to $10$. The real and imaginary parts of the corresponding determinant by the symbolic modified Dodgson＊s algorithm are shown in Figure 3.

\begin{figure}[htbp]
  \centering
  \includegraphics[scale=0.7]{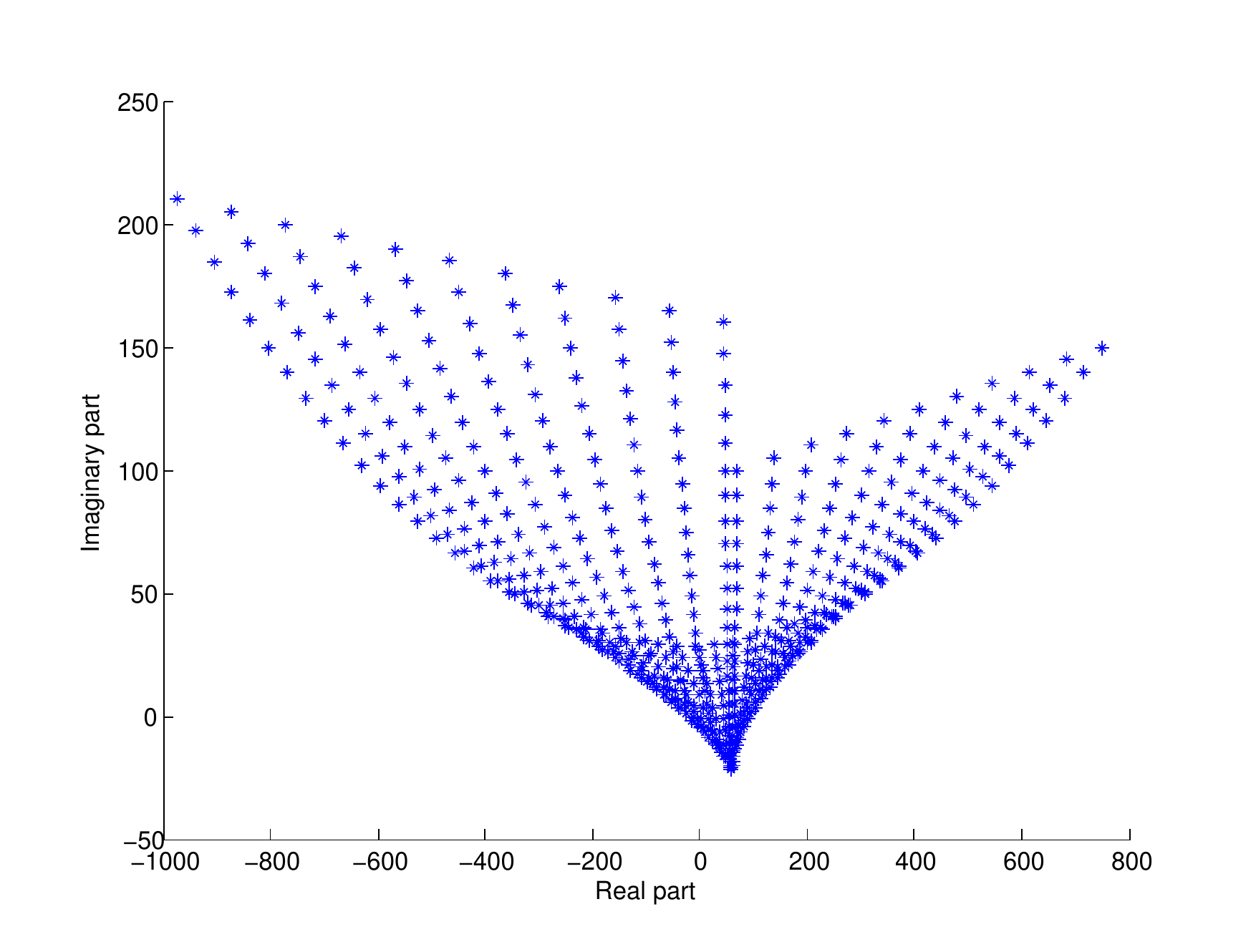}
  \caption{The sensitivity of the symbolic modified Dodgson＊s algorithm to matrix $A_7$.}
  \label{figure:3}
\end{figure}

However, there exist also some matrices which we may replace the zero entries in the interior of intermediate matrix $A^{(n-k)}$ with symbolic variables instead of changing the determinant of original matrix by Dodgson's method. For example,
\[
A = \left[ {\begin{array}{*{20}{c}}
1&0&1&0&1\\
0&1&1&1&1\\
1&2&1&1&2\\
{ - 1}&1&1&2&1\\
0&1&0&1&0
\end{array}} \right].
\]
The special computing process may be not described in detail here. But which of determinants can be calculated like this is still an open problem.

(2) In addition, by the above special examples, we may see that the reason why the symbolic algorithm is successful is that the determinant is
a continuous function on these symbolic variables, which makes us may evaluate its limit when these variables approximate zero.
Therefore, we may also let some non-zero elements of the determinant be some corresponding symbolic variables and
then make these variables approximate their real values to obtain its determinant.

To better understand this method, let us reexamine Example \ref{E2.2}:
\begin{example}\label{E2.22}(\cite{Thesis2011}). Computing the determinant of the following matrix
\begin{equation}\label{2.6}
A_1 = \left[ {\begin{array}{*{20}{c}}
1&{ - 2}&1&2\\
{ 3}&4&4&1\\
6&3&3&4\\
7&5&2&{ - 1}
\end{array}} \right].
\end{equation}

Let
\[\begin{array}{l}
{{\bar A}_1} = \left[ {\begin{array}{*{20}{c}}
1&{ - 2}&1&2\\
3&4&4&1\\
6&x&3&4\\
7&5&2&{ - 1}
\end{array}} \right] \Rightarrow \bar A_1^{(3)} = \left[ {\begin{array}{*{20}{c}}
{10}&{ - 12}&{ - 7}\\
{3x-24}&{12-4x}&{13}\\
{30-7x}&{2x-15}&{-11}
\end{array}} \right]\\
 \Rightarrow \bar A_1^{(2)} = \left[ {\begin{array}{*{20}{c}}
{-x -42}&{-7x -18}\\
{111-22x}&{21+6x}
\end{array}} \right] \Rightarrow \bar A_1^{(1)} = \frac{{-4(40x+93)(x-3)}}{{12-4x}} = 40x+93.
\end{array}
\]

Taking the limit as $x$ approaches $3$, we have
\[
\begin{array}{l}
\det A_1=\mathop {\lim }\limits_{x \to 3} \det ({{\bar A}_1}) = \mathop {\lim }\limits_{x \to 3} (40x+93) = 213.
\end{array}
\]
\end{example}

(3) For multiple zeros case, the implementation of the algorithm using Computer Algebra Systems (CAS) such as MAPLE, MACSYMA, MATHEMATICA, and MATLAB
is also very straightforward and effective. The following two examples illustrate this case.

\begin{example}\label{E2.3} Consider the following fourth order determinant:
\begin{equation}\label{2.9}
|A| = \left| {\begin{array}{*{20}{c}}
1&2&3&4\\
5&0&0&6\\
7&0&0&8\\
9&10&11&12
\end{array}} \right|.
\end{equation}
We use four symbolic parameters $x$, $y$, $z$ and $w$ to avoid handling zeros, that is, we define a new determinant
\begin{equation}\label{2.10}
|\overline{A}| = \left| {\begin{array}{*{20}{c}}
1&2&3&4\\
5&x&y&6\\
7&z&w&8\\
9&10&11&12
\end{array}} \right|.
\end{equation}
The original Dodgson's method yields
\begin{equation}\label{2.11}
\begin{array}{l}
{{\bar A}^{(3)}} = \left[ {\begin{array}{*{20}{c}}
{x - 10}&{2y - 3x}&{18 - 4y}\\
{5z - 7x}&{xw - yz}&{8y - 6w}\\
{70 - 9z}&{11z - 10w}&{12w - 88}
\end{array}} \right] \Rightarrow  \cdots  \Rightarrow \\
{{\bar A}^{(1)}} =  - 24xw + 184x + 128w + 24yz - 136z - 176y + 16.
\end{array}
\end{equation}
Taking the limit as $x,y,z$ and $w$ approach zero, we have
\[\mathop {\lim }\limits_{x \to 0,y \to 0,z \to 0,w \to 0} {{\bar A}^{(1)}} = 16,\]
which shows that the value of the determinant \eqref{2.9} is 16.
\end{example}

Finally, when the zero appears in a $3\times 3$ block of zeros, we know that the double-crossing method may fail.
However, in this case, the symbolic algorithm still succeeds, but it needs more parameters. For example, the following determinant
will require nine parameters to calculate its value by Dodgson's method.
\[|A| = \left| {\begin{array}{*{20}{c}}
1&3&5&7&9\\
{ - 2}&0&0&0&{ - 4}\\
{ - 6}&0&0&0&{ - 8}\\
{ - 10}&0&0&0&{ - 12}\\
9&7&5&3&1
\end{array}} \right|.\]

Finally, let us summarize the above symbolic modified algorithm as follows:
\begin{enumerate}
  \item Firstly, replace the elements that cause entries in the interior of $n\times n$ matrix $A$ to be zeros with the symbolic variables (e.g., $x,y,z,\ldots$).
  \item Find the $2\times 2$ determinant for every four adjacent terms to form a new $(n-1)\times(n-1)$ matrix $A^{(n-1)}$.
  \item Repeat this step to produce an $(n-2)\times (n-2)$ matrix, and then divide each term by the corresponding entry in the interior of
        the original matrix $A$, to obtain matrix $A^{(n-2)}$.
  \item Continue ``condensing" the matrix down, until a single symbolic expression, i.e., $A^{(1)}$ is obtained. This final limit of this symbolic expression as symbolic variables approach the original values in matrix $A$ will be $\det A$.
\end{enumerate}

All in all, this symbol algorithm can work in all cases and preserve the idea of Dodgson's method, which are also suited for implementation using Computer Algebra Systems (CAS) such as MATLAB, MACSYMA, MAPLE and MATHEMATICA.

\section{Concluding remarks}

From the above discussion, one can see that unique utilization of matrix condensation techniques yield an elegant process that has promise for
parallel computing architectures. In this paper, we mainly summarize some progress on the research of Dodgson's method, and some examples are given to clarify these corresponding problems. Numerical experiments show that the symbolic algorithm may be used to overcome the drawback of Dodgson's method by
the aid of computer algebra systems (CAS).

\bibliographystyle{plain}%{unsrt}%{astron}% {astron}% {amsplain}

\bibliography{CMJ-bib-1}

\end{document}